\documentclass[a4paper,11pt,draft]{article}
\usepackage[english]{babel}
\usepackage{amsmath,amsfonts,amsthm}
\usepackage{amssymb}
\usepackage[cp1251]{inputenc}
\newtheorem{theorem}{\sc\bf Theorem.}
\newtheorem{lemma}{\sc\bf Lemma}
\newenvironment{proof1} {\smallskip\normalfont{\scshape \bf Proof of the theorem 1.}}{\vspace{.2 ex}}
\newenvironment{proof2} {\smallskip\normalfont{\scshape \bf Proof of the theorem 2.}}{\vspace{.2 ex}}
\newenvironment{proof3} {\smallskip\normalfont{\scshape \bf Proof of the theorem 3.}}{\vspace{.2 ex}}
\newcommand{\C}{\bold C}

\begin{document}
\begin{center}
 {\bf On the Global Structure of Hopf Hypersurfaces in a Complex Space form.}
 \end{center}
\begin{center}
 {\bf A.A.Borisenko.}
   \end{center}

\hrule \vspace{0.5cm} {\bf Abstract.} It is known that a tube over a K\"ahler submanifold in a
complex space form is a Hopf hypersurface. In some sense the reverse statement is true: a connected
compact generic immersed $C^{2n-1}$ regular Hopf hypersurface in the complex projective space is a
tube over an irreducible algebraic variety. In the complex hyperbolic space a connected compact
generic immersed $C^{2n-1}$ regular Hopf hypersurface is a geodesic hypersphere. \vspace{0.3cm}
 \hrule
 \vspace{0.5cm}
\begin{center}
{\bf Introduction.}
\end{center}
 A natural class of real hypersurfaces in a complex space form $\overline M (c)$
of constant holomorphic curvature 4c is the class of Hopf \ hypersurfaces. For a unit normal vector
$\xi $ of a hypersurface $M$ the vector $J\xi $ is a tangent vector to $M$, where $J$ is the
complex structure of the complex space form $\overline M (c)$.
 \vspace{0.5cm}

{\bf Definition}{\em A hypersurface $M\subset \overline M (c)$ is
called a Hopf hypersurface if the vector $J\xi $ is a principal
direction at every point of $M$.}

 \vspace{0.5cm}

Y.Maeda \cite{[11]} proved that for Hopf hypersurfaces in the
$n$-dimensional complex projective space $\bf {C} P^n$ the
corresponding principal curvature in the direction $J\xi $ is
constant. It is known that a tube over a K\"ahler submanifold in a
complex projective space is a Hopf hypersurface. T.E. Cecil and
P.J. Ryan studied the local and global structure of Hopf
hypersurfaces with constant rank of the focal map $\Phi _r$.

Let $M$ be an embedded hypersurface of $\overline M (c)$ of the
regularity class $C^2$. Let $NM$ be the normal bundle of $M$ with
projection $p:\ NM\to M$ and let $BM$ be the unit normal bundle.
For $\xi \in NM$ let $F(\xi )$ be the point in $\overline M (c)$
reached by traversing a distance $|\xi|$ along the geodesic in
$\overline M (c)$ originating at $x=p(\xi )$ with the initial
tangent vector $\xi $.

A point $P\in \overline M (c)$ is called a focal point of multiplicity $\nu >0$ of $(M,\, x)$ if
$P=F(\xi )$ and the Jacobian of the map $F$ has nullity $\nu $ at $\xi $. \vspace {0.5cm}

 {\bf Definition}{\em The tube of radius $r$ over $M$ is the image of the
map $\Phi _r:\ BM\to \overline M (c)$ given by $\Phi _r(\xi
)=F(r\xi )$, \ $\xi \in BM$.}

\vspace {0.5cm}

T.E. Cecil and P.J. Ryan had proved the following result:

\begin{lemma}  \label{L1}\cite{[1]} Let $M$ be a connected, orientable Hopf
hypersurface of $\C P^n$ with corresponding constant principal
curvature $\mu =2\cot 2r$. Suppose the map $\Phi _r$ has constant
rank $q$ on $M$. Then $q$ is even and every point $x_0 \in M$ has
a neighbourhood $U$ such that $\Phi _r(U)$ is an embedded complex
$q/2$-dimensional submanifold of $\C P^n$.
\end{lemma}

We remark that, in Lemma 1 and Lemma 13 below,  $C^3$ regularity is enough. From Lemmas 1 and 13 we
obtain that Hopf hypersurface with  $\Phi _r$ of constant rank is an analytical hypersurface. It
follows from this fact that $\Phi _r(U)$ is a complex submanifold and parametrizations functions of
$\Phi _r(U)$ satisfy an elliptic system of the $PDE's$ with analytical coefficients. From $C^2$
regularity of $\Phi _r(U)$ we obtain that $\Phi _r(U)$ is analytic.

The global version of Lemma 1 has the following form \cite{[1]}:

Let $M$ be a connected compact embedded real Hopf hypersurface in
$\C P^n$ with corresponding constant principal curvature $\mu
=2\cot 2r$. Suppose the map $\Phi _r$ has constant rank $q$ on
$M$. Then $\Phi _r$ factors through a holomorphic immersion of the
complex $q/2$-dimensional manifold $M/T_0$ into $\C P^n$, where
$T_0$ are $(2n-q-1)$-dimensional spheres, the leaves of the
distribution
$$
T_0(x)=\left\{ y\in T_xM,\ (\Phi _r)_{\ast }(y)=0\right\}.
$$

{\em Acknowledgements.}The author began to study Hopf hypersurfaces when he was a visiting
professor at the University of Valencia (Spain) from January 1997 to June 1997. The author would
like to thank Vicente Miquel  and Yury Nikolaevsky for closely reading a preliminary version of
this paper and for their remarks. The author would like to thank the referee for  valuable remarks.

\begin{center}
{\bf 1. The main results}
\end{center}

The following theorem gives a complete description of the global
structure of Hopf hypersurfaces in complex space forms.

Let $M$ be an immersed regular hypersurface in a regular manifold
$N$. Suppose that for a point $P\in N$ of self-intersection the
linear span of the tangent hyperplanes to the branches of $M$
coincides with tangent space $T_PN$ of the ambient manifold. This
point is called a generic point of self-intersection. If every
point of self-intersection of the hypersurface $M$ is a generic
point of self-intersection then the hypersurface $M$ is called a
generic immersed hypersurface.

\begin{theorem}\label{T1} Let $M$ be a $C^{2n-1}$ regular compact generic
immersed orientable Hopf hypersurface in the complex projective
space ${\C }P^n$ $(n\geqslant 2)$. Then $M$ is a tube over an
irreducible algebraic variety.
\end{theorem}

{\bf Corollary} {\em Let $M$ be a $C^{2n-1}$ regular connected
compact embedded Hopf hypersurface in the complex projective space
$\C P^n \ (n\geqslant 2)$. Then $M$ is a tube over an irreducible
algebraic variety.} \vspace{0.5cm}

The following are some standard examples of Hopf hypersurfaces in
$\C P^n$ of constant holomorphic curvature 4.
\begin{enumerate}
\item A geodesic hypersphere $M$ is the set of points at a fixed
distance $r<\frac{\pi }{2}$ from a point $P\in \C P^n$. It is
obvious that $M$ is also the tube of radius $\frac{\pi }{2} -r$
over the hyperplane $\C P^{n-1}\subset \C P^n$ dual to the point
$P$.

\item A tube over a totally geodesic $\C P^k$ $(1\leqslant k
\leqslant n-1)$.

\item A tube over a totally geodesic real projective space $RP^n$
and over a complex quadric $Q^{n-1}= \{ (z_0,\dots ,z_n\} \subset
\C P^n:\ z^2_0+z^2_1+\dots +z^2_n=0\} $.
\end{enumerate}

A tube of small radius $r$ over a closed irreducible algebraic
manifold in $\C P^n$ is an analytic Hopf hypersurface. But let
$f=x_0^6\, x^2_3+x^3_1\, x^5_2=0$ be the algebraic variety $M$ in
$\C P^3$. The point $P(1,\, 0,\, 0,\, 0)$ is a singular point
($\text{grad}\, f/P=0$). In any neighbourhood of the point $P$ the
normal curvatures at smooth points vary from $-\infty $ to
$+\infty $.  From Lemma 12 below it follows that normal curvatures
of the tube of any radius $r$ tend to $+\infty $. It follows that
the tube of any radius $r$ has regularity less then $C^{1,1}$.

V. Miquel had proved the following theorem:

\vspace{0.5cm}{\bf Theorem}{\em(V. Miquel [13]) Let $M$ be a
connected compact embedded Hopf hypersurface in $\C P^n$ contained
in a geodesic ball of radius $R< \frac{\pi }{2}$.

Suppose that
\begin{enumerate}
\item $M$ has constant mean curvature $H$;

\item The principal curvature $\mu $ in the direction $J\xi $
satisfies the inequality
$$
\mu \geqslant 2\cot \left( 2\text{arc}\,\cot \left[
\frac{(2n-1)H-\mu }{2n-2} \right] \right).
$$

Then $M$ is a geodesic hypersphere.
\end{enumerate}}

We prove the following theorem.

\begin{theorem}\label{T2} Let $M$ be a $C^{2n-1}$ regular connected compact
generic immersed orientable Hopf hypersurface in the complex
projective space $\C P^n$ $(n\geqslant 2)$ contained in a geodesic
ball of radius $R<\frac{\pi }{2}$. Then $M$ is a geodesic
hypersphere.
\end{theorem}

Let $\C H^n$ be the complex hyperbolic space of constant holomorphic curvature $-4$. We prove the
following theorem.

\begin{theorem}\label{T3}  Let $M$ be a connected compact generic immersed
orientable $C^{2n-1}$ regular Hopf hypersurface in the complex
hyperbolic space $\C H^n$ $(n\geqslant 2)$. Then the Hopf
hypersurface $M$ is a geodesic hypersphere.
\end{theorem}

\begin{center}
{\bf 2. Lemmas}
\end{center}

\begin{lemma}\label{L2} (Y. Maeda, \cite{[11]}) Let $M$ be a connected Hopf
hypersurface in the complex projective space $\C P^n$. Then the
principal curvature $\mu $ of $M$ in the direction $J\xi$ is
constant.
\end{lemma}

Let $A_{\xi }$ be the shape operator of $M$.

\begin{lemma}\label{L3}  (T.E. Cecil, P.J. Ryan \cite{[1]}) Suppose $J\xi $ is an
eigenvector of $A_{\xi }$ with an eigenvalue $\mu $. Then we have:

a) $(F_{\ast })_{r\xi }(X,\, 0)=0$ if $\lambda =\text{cot}\, r$ is
an eigenvalue of $A_{\xi }$ and $X$ is a vector in the eigenspace
$T_{\lambda }$ corresponding to the eigenvalue $\lambda $.

b) $(F_{\ast })_{r\xi }(J\xi ,\, 0)=0$ if $\mu =2\text{cot}\, 2r$.

c) $(F_{\ast })_{r\xi }(X,\, V)\ne 0$ except as determined by (a)
and (b).
\end{lemma}

Now, let $M$ be a real hypersurface of a complex space form
$\overline M ^n(c)$ of constant holomorphic curvature $4c$  and
let $\xi $ be a unit normal field on $M$. If $X\in T_PM$, $P\in
M$, then one has a decomposition
$$
JX=\phi X+f(X)\xi
$$
into the tangent and normal components respectively. So, $\phi $
is a $(1,\, 1)$-tensor field and $f$ is a $1$-form. Then they
satisfy
$$
\phi ^2X=-X+f(X)U, \quad \phi U=0, \quad f(\phi X)=0
$$
for any vector field $X$ tangent to $M$, where $U=-J\xi $.
Moreover, we have
$$
g(\phi X,\, Y)+g(X,\, \phi Y)=0, \quad f(X)=g(X,\, U);
$$
$$
g(\phi X,\, \phi Y)=g(X,\, Y)-f(X)f(Y)
$$
with $g$ the metric tensor in $\overline M ^n(c)$. We denote by
$A$ the shape operator on $T_PM$ associated with $\xi $.

\begin{lemma}\label{L4}

1.(\cite{[9]}) Let $M$ be a Hopf hypersurface in $\overline
M^n(c)$. Then we have
\begin{enumerate}
\item[a)]\quad $-2c\phi =\mu (\phi A+A\phi )-2A\phi A; $
\item[b)]\quad $ X\mu =(U\mu )f(X)$
\end{enumerate}
 and

$$
(U\mu )\, g( (\phi A+A\phi )X,\, Y) =0,
$$
where $\mu $ is the principal curvature in the direction $U=-J\xi
$, $X,\ Y$ are vectors tangent to $M$, and $U\mu $ is the
derivative of the function $\mu $ in the direction $U$. Moreover,
if $\phi A+A\phi =0$  then
$$
cg(X,\, \phi Y)=-g(\phi AX,\, AY)=g(A\phi X,\, AY),
$$
$$
cg(\phi X,\, \phi X)=-g(A\phi X,\, A\phi X)
$$
and so $c\leqslant 0$.

2.(\cite{[11]}) Let $M$ be a Hopf hypersurface in $\C P^n$. If $X\in T_{\alpha }\subset T_PM$, then
$$
JX \in T_{\mu \alpha +2/2\alpha -\mu }\subset T_PM,
$$
where $T_{\alpha}$ is an eigenspace corresponding to a principal
curvature $\alpha $.
\end{lemma}

It follows from the equation (a) of the first part of the lemma that $\alpha$ cannot be equal to
$\mu$ or to $\mu/2$.

{\bf Definition} {\em Let $A$ be a subset of a metric space $X$.
Let $\delta (A)$ denote the diameter of $A$, and let
$$
\delta ^p(A)=[\delta (A)]^p\ \text{for }p>0,$$
$$
\delta ^0(A)= \left\{\begin{array}{cc}
 1, & \text{if } A\ne \varnothing ;\\
             0  & \text{if } A=\varnothing .
\end{array}\right.
$$

For $p\geqslant 0$ and $\varepsilon >0$ define.
$$
\aligned &H^p_{\varepsilon }(A)=\inf \left\{ \sum ^{\infty }_{i=1}
\delta ^p(A_n):\
     A\subset \cup A_n\ \text{and } \delta (A_n)<\varepsilon \right\} ;\\
&H^p(A)=\lim _{\varepsilon \to 0^{+}} H^p_{\varepsilon
}(A)=\text{sup} \,
    H^p_{\varepsilon }(A).
\endaligned
$$
We call $H^p$ the Hausdorff $p$-measure.}

\begin{lemma} \label{L5} (H. Federer, \cite{[4]}) If $m>\nu \geqslant 0$ and
$k\geqslant 1$ are integers, $A$ is an open subset of $R^m$,
$B\subset A$, $Y$ is a normed vector space and $f:\ A\to Y$ is a
map of class $C^k$ such that
$$
\text{Dim im} \, f_{\ast }(x)\leqslant \nu \quad \text{for} \ x\in
B,
$$
then
$$
H^{\nu +(m-\nu)/k}[f(B)]=0.
$$
\end{lemma}

{\bf Definition} {\em Let $\Omega $ be a complex manifold. A set
$A\subset \Omega $ is called an analytic set in $\Omega $ if for
each point $a\in \Omega$, there exist a neighbourhood $U$ of $a$
and functions $f_1,\, \dots ,\, f_N$ holomorphic in $U$ such that
$A\cap U=Z_{f_1}\cap \dots \cap Z_{f_k}\cap U$, where $Z_f$ is the
set of zeros of a holomorphic function~$f$.}

\vspace{0.5cm}

A point $a$ of an analytic set $A$ is called a regular point if there exists a neighbourhood $U$ of
$a$ in $\Omega $ such that $A\cap U$ is a complex submanifold of $U$. The complex dimension of
$A\cap U$ is then called the dimension of $A$ at the point $a$ and is denoted by $\dim _{a}A$. The
set of all regular points of an analytic set is denoted by $\text{reg}\, A$. Its complement
$A\setminus \text{reg}\, A$ is denoted by $\text{sng}\, A$. The set $\text{sng}\, A$ is called the
set of singular points of the set $A$. It can be shown by induction on the dimension of the
manifold $\Omega $ that $\text{sng}\, A$  is nowhere dense and closed. This allows us to define the
dimension of $A$ at any point $a$ of $A$ as
$$
\dim _a A=\lim_{z\to a}\, \dim _z A\ (z\in \text{reg}\, A).
$$

The set $A$ is called purely $p$-dimensional if $\dim _z A=p$ for
all $z\in A$ \cite{[2]}, \cite{[3]}.

\begin{lemma} \label{L6} (B. Shiffman, \cite{[16]}) Let $E$ be a closed subset of a
complex manifold $\Omega $ and let $A$ be a purely $q$-dimensional
analytic subset of $\Omega \setminus E$. If $H^{2q-1}(E)=0$ then
the closure $\overline A$ of the set $A$ in $\Omega $ is a purely
$q$-dimensional analytic subset of $\Omega $.
\end{lemma}

{\bf Definition}{\em (D.Mumford, \cite{[14]}) Let $U\subset \C ^n$
be an open set. A closed subset $X\subset U$ is a $\ast $-analytic
subset of $U$ if $X$ can be decomposed
$$
X=X^{(r)}\cup X^{(r-1)}\cup \dots \cup X^{(0)},
$$
where for all $i$, $X^{(i)}$ is an $i$-dimensional complex
submanifold of $U$ and $\overline X ^{(i)}\subset X^{(i)}\cup
X^{(i-1)}\dots \cup X^{(0)}$. If $X^{(r)}\ne \varnothing $, then
$r$ is called the dimension of $X$.}

An analytic set is always $\ast$-analytic \cite{[14]}.

\begin{lemma} \label{L7} (Chow's Theorem, \cite{[14]}) If $X\subset \C P^n$ is a
closed $\ast$-analytic subset, then $X$ is a finite union of
algebraic varieties.
\end{lemma}

\begin{lemma} \label{L8} \cite{[3]} An analytic set $A$ in a complex manifold $\Sigma
$ is irreducible if and only if the set $\text{reg}\, A$ is
connected.
\end{lemma}

Let $X\subset \C P^n$ denote a closed irreducible algebraic variety of dimension $l$ which may have
singularities and let $X_e\subset X$ denote the non-empty open subset of its smooth points. For the
definitions of irreducible singular and smooth points see \cite{[14]}. Define
$$
V'_{X}=\left\{ (x,\, y)\in \C P^n\times  \C \breve P^n\, |\, x\in
X_e\ \text{and }y\text{ is tangent hyperplane at } x \right\} ,
$$
where $\C \breve P^n$ is the dual complex projective space.

The closure $V_X$ of $V'_X$ on Zariski topology in $\C P^n\times
\C \breve P^n$ is called the tangent hyperplane bundle of $X$.
 It is a closed irreducible algebraic variety of dimension
$(n-1)$. The first projection maps $V_X$ onto $X$
$$
\pi _1\colon\ V_X \to X,\quad (x,\, y)\to x.
$$
Consider now the second projection
$$
\pi _2:\ V_X \to \C \breve P^n,\quad (x,\, y)\to y.
$$
Its image $\breve X=\pi _2(V_X)$ is a closed irreducible variety
of $\C \breve P^n$ of dimension at most $(n-1)$, the dual variety
of $X$ \cite{[9]}.

\begin{lemma} \label{L9} (Duality Theorem) \cite{[6]}, \cite{[10]} The tangent hyperplane
bundles of an closed irreducible algebraic variety $X$ and its dual variety $\breve X$ coincide: WE
have $ V_{\breve X}=V_{X}\ \text{and hence }\ \Breve {\Breve X} =X. $
\end{lemma}

Let $\C P^n$ be the complex projective space with standard
Fubini-Study metric. To a hyperplane $L\subset \C P^n$ passing
through a point $x\in \C P^n$ we associate the point $y\in \C P^n$
representing the complex line  in $\C ^{n+1}$ orthogonal to $L$.
Then the distance $\rho (x,y)$ is equal to $\pi /2$. One can
identify $\C \breve P^n$ with $\C P^n$ in this way and consider
$\breve X$ as a subset in $\C P^n$.

 It is possible to define a
tube over an closed irreducible algebraic variety $X\subset \C P^n$ which may have singularities.
Let $(x,\, y)\in V_X\subset \C P^n\times \C \breve P^n =\C P^n\times  \C P^n$, $x\in X$, $y\in
\breve X$, and let $L(x,\, y)$ be a complex projective line through $x,\ y\in \C P^n$. Then
$L(x,y)$ is a totally geodesic two-dimensional sphere in $\C P^n$ of curvature 4, the distance
$\rho (x,\, y)$ is equal to $\pi /2$, and $x$ and $ y$ are poles of the sphere $L(x,\, y)$. The set
of points of $L(x,\, y)$ at a distance $r$ from the point $x$ is a circle $S_r(x,\, y)$ with the
center $x$. The union
$$
S_r=\bigcup _{(x,y)\in V_X} S_r(x,\, y)
$$
is called the tube of radius $r$ over $X$. The set $S_r$ is the
tube of radius $\frac{\pi }{2}-r$ over the dual variety $\breve
X$.

If all the points of $X$ are regular this definition coincides
with one above.

The set of points $\text{sng}\, V_X\subset V_X$ such that $(x,\,
y)\in \text{sng}\, V_X$ if $x\in \text{sng}\, X$ or $y\in
\text{sng}\, \breve X$ is a closed algebraic subvariety of $V_X$,
$\text{reg}\, V_X=V_X\setminus \text{sng}\, V_X$ is an open set of
$V_X$ in the Zariski topology.

Let $X\subset \C P^n$ be a closed irreducible algebraic variety
and let $x_0$ be a Zariski open set in $X$. Then the closure of
$x_0$ in the classical topology is $X$ \cite{[14]}.

Let us take the Segre map
$$
\sigma :\ \C P^n\times \C \breve P^n \to \C P^{(n+1)^2-1}.
$$
Then $\sigma (V_X)$ is a closed irreducible algebraic variety in $\C P^{(n+1)^2-1}$ and the set
$\text{reg}\, V_X$ is an open set of $V_X$ in the Zariski topology.

As corollary we obtain the following result

\begin{lemma} \label{L10} The closure of the set $\text{reg}\, V_X\subset \C
P^n\times  \C P^n$ in the standard topology coincides with the
tangent bundle $V_X$.
\end{lemma}

Therefore the tube over $X$ is the closure of the set
$$
\bigcup _{(x,y)\in \text{Reg}V_X} S_r(x,\, y)
$$

\begin{lemma} \label{L11} \cite{[5]} Let $X$ be a compact topological space. Suppose
$A$ is a closed subset such that $X\setminus A$ is a smooth
$n$-dimensional orientable manifold without boundary. Then
$$
H_q(X,\, A)\backsimeq H^{n-q}(X\setminus A),
$$
where $H_i,\ H^i$ are homology and cohomology groups.
\end{lemma}

\begin{lemma} \label{L12} \cite{[1]} Suppose $J\xi $ is an eigenvector of the shape
operator $A_{\xi }$ of a Hopf hypersurface $M$ in the complex
projective space, with the corresponding eigenvalue $2\cot
2\Theta$, $0<\Theta <\frac{\pi }{2}$. Suppose $J\xi ,\, X_2,\,
\dots ,\, X_n$ is a basis of principal vectors of $A_{\xi }$ with
$A_{\xi }X_j=\text{cot}\, \Theta _jX_j$, $2\leqslant j\leqslant
n$, $0<\Theta _j<\pi $; $\frac{\partial }{\partial t_j}$
$(2\leqslant j \leqslant k)$ are normal vectors. Then the shape
operator $A_r$ of the tube $\Phi _r$ is given in terms of its
principal vectors by

(a) $A_r\left( \frac{\partial }{\partial t_j}\right) =
     -\cot \, r\left(\frac{\partial }{\partial t_j}\right) ,\ \
      2\leqslant j\leqslant k;$

(b) $A_r\left( X_j,\, 0\right) =\cot \left( \Theta _j-r\right)
     \left( X_j,\, 0\right) ,\ \
      2\leqslant j\leqslant n;$

(c) $A_r(J\xi ,\, 0)=\cot \left( 2(\Theta -r)\right) (J\xi ,\,
0)$.
\end{lemma}

For a complex hyperbolic space $\C H^n$ the following analog of
 Lemma 1 holds:

\begin{lemma} \label{L13}\cite{[13]} Let $M$ be an orientable Hopf hypersurface of
$\C H^n$ such that the principal curvature $\mu $ in the direction $J\xi $ is constant and equal to
$\mu =2\coth 2r$. Suppose that $\Phi _r$ has constant rank $q$ on $M$. Then for every point $x_0\in
M$ there exists an open neighbourhood $U$ of $x_0$ such that $\Phi _rU$ is a $q/2$-dimensional
complex submanifold embedded in $\C H^n$.
\end{lemma}

\begin{lemma}\label{L14}  \cite{[15]} Let $\Omega $ be a Hermitian complex manifold
with exact fundamental form $\omega=d\gamma $. Let $A$ be an
analytical $q$-dimensional set with boundary $\partial A\subset
\Omega $ such that $A\cup \partial A$ is compact.

Then
$$
H^{2q}(A)\leqslant \frac1{q}\left( \text{max} _{\partial A}|\gamma
|\right) \, H^{2q-1}(\partial A),
$$
where $H^{2q}(A)$, $H^{2q-1}(\partial A)$ are Hausdorff measures, and

$$|\gamma |(z)=\max \left\{ |\gamma (\upsilon )|:\, \upsilon \in
      T_z\Omega ,\ |\upsilon |=1\right\}. $$
\end{lemma}

\begin{lemma} \label{L15} \cite{[8]} Let $M$ be a Hopf hypersurface of a complex
space form $\overline M ^n(c)$ $(c\ne 0)$. If $U$ is an
eigenvector of $A$, then the principal curvature $\mu =g(AU,\, U)$
is constant.
\end{lemma}

\begin{center}
{\bf 3. Proofs of the Theorems}
\end{center}

Let $M_s$ be the set of points of $M$ such that $\text{rank}\, (\Phi _r)_{\ast }(M_s)=s$, $F_s=\Phi
_r(M_s)$, $F=\Phi _r(M)$. From Lemma 4 we obtain that if $X\in T_{\alpha }\subset T_PM$ where
$T_{\alpha }$ is the eigenspace corresponding to the principal curvature $\alpha =\text{cot}\, r$,
then $JX\in T_{\alpha }$. Hence $s$ is even and if $s<2q$, then $s\leqslant 2q-2$.

Let
$$
E=\bigcup _{s<2q} F_s\cup F_0
$$
$$
F_0=\left\{ x\in F\colon x=\Phi _r(L_1)=\Phi _r(L_2),\, L_1\ne
L_2\subset M, \text{rank} \, (\Phi _r)_{\ast }(P_1)=\text{rank} \,
(\Phi _r)_{\ast }(P_2)=2q \right\},
$$
 for $ P_i \in L_i$,  where $ L_i $ are leaves of the
distribution $ Ker (\Phi_r )_*$.

\begin{proof1}
  Let $M$ be a compact Hopf
hypersurface in $\C P^n$. This means that the vector $J\xi $ is a principal direction of $M$, where
$\xi $ is the unit normal vector and $J$ is the complex structure in $\C P^n$. From Lemma 2 it
follows that the corresponding principal curvature $\mu $ is constant, $\mu =2\cot 2r$. Let $2q$ be
the maximal rank of $(\Phi _r)_{\ast }$ on $M$. Let $P\in M$ be a point such that $\text{rank}\,
(\Phi _r)_{\ast }(P)=2q$ and let $M_{2q}$ be the corresponding connected component of $M$ such that
$P\in M_{2q}$ and for $Q\in M_{2q}$ \ $\text{rank}\, (\Phi )_{\ast }(Q)=2q$. Set $F_{2q}=\Phi
_r(M_{2q})$, $\widetilde F =F_{2q}\cap (\C P^n\setminus E)$. From Lemma \ref{L1} we obtain that
$\widetilde F$ is a purely analytic set, $\text{dim} _z\widetilde F =q$, $z\in \widetilde F$.

Locally $F_0$ is a transversal intersection of two complex
submanifolds of dimension $q$. Hence $F_0$ is an analytic set of
real dimension $\leqslant 2q-2$. Then its Hausdorff measure
$$
H^{2q-1}(F_0)=0.
$$

Now apply Lemma \ref{L5} to the set $E_1=\bigcup\limits _{s<2q}
F_{s}$ and the map $\Phi _r$. Then $\nu \leqslant 2q-2$.

If the class of regularity of $M$ is greater or equal to
$2(n-q+1)$ then the class of regularity of $\Phi _r$ is
$k\geqslant 2(n-q+1)-1$ and
$$
\nu +\frac{2n-1-\nu }{k} \leqslant 2q-2+\frac{2n-1}{k} \leqslant
2q-1,
$$
for $k\geqslant 2n-1$. From Lemma \ref{L5} we have
$H^{2q-1}(E_1)=0$ and so $H^{2q-1}(E)=0$. From Lemma \ref{L6} we
obtain that the closure of $\widetilde F$ is a purely
$q$-dimensional analytic subset of $\C P^n$. Since any analytic
subset is $\ast $-analytic we get from Chow's Theorem (Lemma
\ref{L7}) that $\text{cl}\, \widetilde F \subset \C P^n$ is a
finite union of algebraic varieties. An analytic set $A$ is an
irreducible if and only if the set $\text{reg}\, A$ is connected.
From Lemma 8 it follows that $\text{cl}\, \widetilde F$ is
irreducible as analytic set and we obtain that $\text{cl}\,
\widetilde F =X$ is an irreducible algebraic variety.

Let $S_r$ be a tube over $X=\text{cl}\, \widetilde F$. From Lemma
10 we have $S_r\subset M$ and $S_r=\text{cl}\, M_{2q}$. We will
prove that $\text{cl}\, M_{2q}=M$. Suppose that $\text{cl}\,
M_{2q}\ne M$. Then in every neighbourhood of a point $P\in
\partial M_{2q}$ there exist points $Q\in M\setminus \text{cl}\,
M_{2q}$. Let $P\in \partial M_{2q}$. Then $P\in S_r(x,\, y)$ such
that $x\in \text{sng}\, X$, $y\in \text{sng}\, \breve X$. Then
$$
\partial \, M_{2q}=\bigcup\limits_{ x\in\, \text{sng} \, X,
   y\in \,\text{sng} \, \breve X } S_r(x,\, y).
$$
Otherwise some neigbourhood of $P$ belongs to $\text{cl}\, M_{2q}$
and $P\in \text{int}\, \text{cl}\, M_{2q}$. The set of points
$$
\text{sng}\, (X,\, \breve X)=\text{sng}\, X\times \C P^n \cap \C
P^n\times \text{sng}\, \breve X \subset V_X\subset \C P^n \times
\C P^n
$$
is a closed algebraic subvariety of $V_X$. The dimension of $\text{sng}\, (X,\, \breve X)\leqslant
n-2$ because the dimension of $V_X$ is equal to $n-1$. The set $\partial \, M_{2q}$ is a fiber
bundle over $\text{sng}\, (X,\, \breve X)$ with the circle $S^1$ as a leaf. The real dimension of
$\text{sng}\, (X,\, \breve X)$ is  $\leqslant 2(n-2)$ whence
$$
H_{2n-3}\left( \text{sng}\, (X,\, \breve X),\, \bold Z\right) =0.
$$
For $E=\partial \, M_{2q}$, $B=\text{sng}\, (X,\, \breve X)$,
$F=S^1$ the exact Thom-Gysin sequence has the form \cite{[17]}
$$
\aligned &H_{2n-1}\left( \text{sng}\, (X,\, \breve X),\, \bold
Z\right) \to
H_{2n-3}\left( \text{sng}\, (X,\, \breve X),\, \bold Z\right) \to \\
&\to H_{2n-2}\left( \partial \, M_{2q},\ \bold Z \right) \to
H_{2n-2}\left( \text{{sng}}\, (X,\, \breve X),\, \bold Z\right) ,
\endaligned
$$
$$
0\to 0\to H_{2n-2}\left( \partial \, M_{2q},\ \bold Z \right) \to
0.
$$
We obtain
$$
H_{2n-2}\left( \partial \, M_{2q},\ \bold Z \right) =0.
$$

Next,we apply Lemma \ref{L11} with $X=M$, $A=\partial \, M_{2q}$. Then
$$
H_{2n-1}\left(M,\ \partial \, M_{2q} \right) = H^0\left(
M\setminus \partial \, M_{2q} \right) .
$$
But $M\setminus \partial M_{2q}$ has $m>1$ connected components
and
$$
H^0\left( M\setminus \partial \ M_{2q},\, \bold Z \right) =
\bigoplus ^{m}_{i=1} \bold Z
$$
is the direct sum of $m$ copies of $\bold Z$ [17].

For the pair $\left(M,\ \partial \, M_{2q} \right) $ the exact
homology sequence has the following form
$$
\aligned &  H_{2n-1}\left(\partial \, M_{2q},\, \bold Z \right)
\to
  H_{2n-1}\left(M,\, \bold Z \right) \to
  H_{2n-1}\left(M,\ \partial \, M_{2q},\, \bold Z \right) \to \\
& \to H_{2n-2}\left(\partial \, M_{2q},\ \bold Z \right) ;
\endaligned
$$
$$
 H_{2n-1}\left(\partial \, M_{2q},\, \bold Z \right) =
  H_{2n-2}\left(\partial \, M_{2q},\, \bold Z \right) =0; \quad
  H_{2n-1}\left(M,\, \bold Z \right) =\bold Z.
$$
It follows that $H_{2n-1}\left(M,\ \partial \, M_{2q},\, \bold Z \right) =\bold Z$. This
contradicts to above result. Thus $\text{cl}\, M_{2q}=M$ and $M$ is a tube over the irreducible
algebraic variety $\text{cl}\, \widetilde F=X$.
\end{proof1}

\begin{proof2}
 Let $S$ be the hypersphere of the
minimal radius $r_0$ such that the hypersurface $M$ is contained in the ball $D$ with  boundary
$\partial D\, =S$. Let $P$ be a point of tangency of $M$ and $S$. Let $\xi $ be the inward unit
normal vector at the point $P$. Then the principal curvature in the direction $J\xi $ is $\mu =
2\text{cot} \, 2\rho \geqslant 2\text{cot} \, 2r_0$, and so $\rho \leqslant r_0<\pi /2$. Another
principal curvature $k_i=\text{cot} \, \Theta _i$ at the point $P$ satisfies the conditions
$\text{cot} \, \Theta _i\geqslant \text{cot} \, r_0$, where $2\text{cot} \, 2r_0$, $\text{cot} \,
r_0$ are principal curvatures of the hypersphere $S$. Then $\Theta _i\leqslant r_0$. Let $r=\rho
-\pi /2$. From Lemma \ref{L12} we obtain that the principal curvatures of the tube $\Phi _r$ over
$M$ are equal to
$$
(k_i)_r=\text{tg} \, (\rho -\Theta _i)\leqslant \text{tg} \,
(r_0-\Theta _i)< \infty .
$$
Hence $\text{rank} \, (\Phi _r)_{\ast }(P)=2(n-1)$ and from Theorem 1 we get that $\Phi
_r(M)=\text{cl}\, \widetilde F =X$ is an irreducible hypersurface of degree $d$. Let $X_k$ be a
sequence of smooth algebraic hypersurfaces such that $\lim X_k=X$, degree $X_k=d$ \cite{[7]}, and
let $\breve X$, $\breve X _k$ be dual algebraic varieties. Then
$$
M=\Phi _{\frac{\pi }{2}-r}(X)=\Phi _r(\breve X)
$$
and from Lemma \ref{L9} we get that  $\breve X =\lim \breve X _k$.
From the above for $\Phi _{\frac{\pi }{2}-r}(X_k)=M_k$,
$$
\lim M_k =M.
$$
For large $k$, $M_k$ is contained in the balls $D_k$ of radius
$R<\pi /2$ and $M_k$ does not intersect complex projective space
$x_0=0$.

Let $f=0$ be the equation of the algebraic hypersurface $X_n$ where $f$ is a homogeneous
polynomial, $\text{grad}\, f\ne 0$. By Bezou Theorem \cite{[15]} the system of equations
$$
x_0=0, \quad f=0, \quad f_{x_0}=0
$$
has a nontrivial solution for $n$ is $\geqslant 3$ and degree of the polynomial $f\geqslant 2$.
This means that $M_k$ intersects the hyperplane $x_0=0$. It follows that $f$ is a linear function
and the $X_k$ are hyperplanes, $M_k$ are hyperspheres. Then the hypersurface $M$ is a geodesic
hypersphere too.

For $n=2$ the equation of the tube has the following parametric
form
$$
z_j=x_j\, \cos r+\sin r \, \frac{\frac{\overline {\partial
f}}{\partial x_j}} {| \text{grad }f|}\, e^{it};
$$
$x_j$ are coordinates of points of the algebraic variety, $0\leqslant t\leqslant 2\pi $;  $0
\leqslant r\leqslant \frac{\pi }2$, $r$ is radius of the tube $\Phi _r$; $j=0,\, 1,\, 2$.

From the real point of view $X$ is a compact two-dimensional
manifold.

Denote
$$
g_1=|x_0\, \cos r|, \quad g_2= \left| \frac{\frac{\overline
{\partial f}}{\partial x_0}} {| \text{grad }f|}\, \sin r\right|,
$$

If the degree of the polynomial $f$ is  $\geqslant2$ the zero sets of these regular functions on
the manifold $X$ are non empty on the manifold $X$. Hence there exists a point $P\in X$ such that
$g_1=g_2=\rho $. Then $z_0=\rho \, \left( e^{i\alpha }+e^{i(\beta +t)}\right) $. Moreover, if
$t=\alpha -\beta -\pi $ then $z_0=0$.

This means that $M_k$ intersects the hyperplane $x_0=0$.

Thus $f$ is a linear function and $M_k$ and $M$ are  geodesic
hyperspheres as in the case $n\geqslant 3$.
\end{proof2}

\begin{proof3}
 Let $S$ be the hypersphere of the
minimal radius $r_0$ such that the hypersurface $M$ is contained in the ball $D$ with boundary $S$.
Let $P_0$ be a point of tangency of $M$ and $S$. Let $\xi $ be the inward unit normal vector of $M$
at the point $P_0$. From Lemma \ref{L15} it follows that the principal curvature $\mu $ in the
direction $J\xi $ is constant. At the point $P_0$ this curvature satisfies the inequality $\mu
\geqslant 2\coth 2r_0$ and $\mu = 2\coth 2r$. We now follow the proof of Theorem 1, using Lemma
\ref{L13} instead of \linebreak Lemma \ref{L1}. Consider the map $\Phi _r$. For a Hopf hypersurface
$\text{rank} \, (\Phi _r)_{\ast }$ is always even. This follows from Lemma \ref{L4}.

Suppose $2q$ is the maximal rank of $(\Phi _r)_{\ast }$ at the points of $M$. Let $P\in M$ be a
point such that $\text{rank} \, (\Phi _r)_{\ast }(P)=2q$ and $M_{2q}$ is the connected component of
$M$ such that for $Q\in M_{2q}$ \ $\text{rank} \, (\Phi _r)_{\ast }(Q)=2q$. As in the proof of
Theorem \ref{T1}, set
$$
F=\Phi _r(M), \quad F_{2q}=\Phi _r(M_{2q}), \quad F_s=\Phi
_r(M_s),
$$
$$
E=F_0\bigcup _{s<2q} F_s; \quad \widetilde F =F_{2q}\cap \C
H^n\setminus E.
$$
We obtain that $\text{cl} \, \widetilde F =X$ is a compact analytic set in $\C H^n$ with boundary
$\partial X\subset E$. The Hausdorff measure $H^{2q-1}(\partial X)=0$. From Lemma \ref{L14} it
follows that $H^{2q}(X)$ is equal to $0$. This is possible only if $q=0$ and $X$ is a point. Then
$M$ is a tube over a point and $M$ is a geodesic hypersphere. \end{proof3}
\begin{center}
\renewcommand{\refname}{\bf {References.}}
  
 \end{center}
\end{document}